\documentclass[a4paper,11pt]{amsart}
\usepackage{amssymb}

\title{Inversion of adjunction on log canonicity}
\author{Masayuki Kawakita}

\address{Research Institute for Mathematical Sciences, Kyoto University,
Kyoto 606-8502, Japan}
\email{masayuki@kurims.kyoto-u.ac.jp}

\newtheorem*{theorem}{Theorem}
\newtheorem{lemma}{Lemma}

\newcommand{\bR}{{\mathbb R}}
\newcommand{\cO}{{\mathcal O}}
\newcommand{\fa}{{\mathfrak a}}
\newcommand{\fm}{{\mathfrak m}}

\newcommand{\rd}[1]{\lfloor{#1}\rfloor}
\newcommand{\ru}[1]{\lceil{#1}\rceil}

\begin{document}
\maketitle

\begin{abstract}
We prove inversion of adjunction on log canonicity.
\end{abstract}

We work over an algebraically closed field of characteristic zero.
Let $X$ be a normal variety and $B$ an effective $\bR$-divisor on $X$.
We call the pair $(X, B)$ a \textit{log pair}
if $K_X+B$ is an $\bR$-Cartier $\bR$-divisor.
A valuation $v$ of the function field of $X$ is called
an \textit{algebraic valuation} of $X$
if it is defined by a prime divisor $E$ on a non-singular variety $\bar{X}$
equipped with a birational morphism $f \colon \bar{X} \to X$.
Such a valuation is denoted by $v_E$.
If $(X, B)$ is a log pair,
we can write $K_{\bar{X}}+\bar{B} = f^*(K_X+B)+A$,
where $\bar{B}$ is the birational transform of $B$ on $\bar{X}$
and $A$ is an exceptional $\bR$-divisor.
We define the \textit{discrepancy} $a_{K_X+B}(v)$ of $v$
as $\textrm{the coefficient of $E$ in $A-\bar{B}$}$.
We say that $(X, B)$ is \textit{log canonical} if
$a_{K_X+B}(v) \ge -1$ for any algebraic valuation $v$.

Let $S$ be a reduced divisor on $X$
which has no common component with the support of $B$,
and let $\nu \colon S^\nu \to S$ denote the normalisation of $S$.
If $(X, S+B)$ is a log pair,
we can construct an effective $\bR$-divisor $B^\nu$ on $S^\nu$,
called the \textit{different} of $B$ on $S^\nu$, in the following manner;
see \cite[Chapter 16]{K+92} or \cite[\S 3]{S92} for details.
The construction is reduced to the case where $X$ is a surface.
Take a resolution of singularities $f \colon \bar{X} \to X$ of the surface $X$
such that the birational transform $S^\nu$ of $S$ on $\bar{X}$
is also non-singular.
Then $S^\nu$ is the normalisation of $S$
and we can find an exceptional $\bR$-divisor $A$ on $\bar{X}$
such that $K_{\bar{X}}+S^\nu+\bar{B} \equiv_f A$.
The different $B^\nu$ is defined as the $\bR$-divisor $(\bar{B}-A) |_{S^\nu}$.
The pair $(S^\nu, B^\nu)$ is a log pair and
it satisfies that $K_{S^\nu}+B^\nu=\nu^*((K_X+S+B) |_S)$.

In this paper we prove \textit{inversion of adjunction}
on log canonicity as stated below.
It was conjectured and proved in dimension three by Shokurov in \cite{S92}.

\begin{theorem}
Let $(X, S+B)$ be a log pair such that
$S$ is a reduced divisor
which has no common component with the support of $B$,
let $S^\nu$ denote the normalisation of $S$,
and let $B^\nu$ denote the different of $B$ on $S^\nu$.
Then $(X, S+B)$ is log canonical near $S$ if and only if
$(S^\nu, B^\nu)$ is log canonical.
\end{theorem}

Note that this statement has been conjecturally generalised to
the precise version of inversion of adjunction in \cite[17.3 Conjecture]{K+92}.
The exposition \cite[Section 5.4]{KM98} is a good reference on
inversion of adjunction.

We begin to prove the theorem.
The statement is local so
we discuss over the germ at a closed point $x \in S \subset X$.
For a resolution of singularities $f \colon \bar{X} \to X$ which is proper,
we call $f$ a \textit{good resolution} if
its exceptional locus on $\bar{X}$ is a divisor
of which the summation with the support of the birational transform of $S+B$
is a simple normal crossing divisor.
Given an ideal sheaf $\fa$ on $X$,
we say that $f$ is a \textit{good resolution for} $\fa$
if furthermore the ideal sheaf $\fa\cO_{\bar{X}}$ on $\bar{X}$ is
an invertible sheaf.
The existence of good resolutions is guaranteed by
Hironaka's resolution theorem \cite{H64}.
We fix a good resolution of singularities $f \colon \bar{X} \to X$.
Let $\bar{S}+\bar{B}$ denote the birational transform of $S+B$ on $\bar{X}$,
$\nu \colon S^\nu \to S$ the normalisation morphism,
and $g \colon \bar{S} \to S^\nu$ the induced morphism.
Let $\{E_i\}_{i \in I}$ denote the set of prime divisors on $\bar{X}$ which
appears in either the exceptional locus of $f$ or the support of $\bar{B}$,
and write $K_{\bar{X}}+\bar{S} = f^*(K_X+S+B)+\sum_{i \in I}a_iE_i$
with $a_i:=a_{K_X+S+B}(v_{E_i})$.
Then $K_{\bar{S}} = g^*(K_{S^\nu}+B^\nu)+\sum_{i \in I_S}a_iE_i|_{\bar{S}}$,
where $I_S := \{i \in I \mid E_i \cap \bar{S} \neq \emptyset\}$.
We mention that $(X, S+B)$ is log canonical if and only if
$a_i \ge -1$ for any $i \in I$,
and $(S^\nu, B^\nu)$ is log canonical if and only if
$a_i \ge -1$ for any $i \in I_S$.
In particular the only-if-part of the theorem is trivial.
It suffices to prove the if-part of the theorem.

From now on we suppose that $a_i \ge -1$ for any $i \in I_S$.
Let $L^\nu$ denote the reduced locus in $S^\nu$ defined as
the union of $g(E_i|_{\bar{S}})$ for all $i \in I_S$ with $a_i=-1$.
$L^\nu$ is called the \textit{locus of log canonical singularities}
for $(S^\nu, B^\nu)$ and it is independent of the choice of $f$.
We let $L$ denote the reduced locus $\nu(L^\nu)$
and $\fa_L$ the ideal sheaf on $X$ defining $L$.
We construct a sequence of ideal sheaves
$\{\fa_j\}_{j\ge0}$ on $X$ in the following manner.

\begin{enumerate}
\item
Set $\fa_0:=\fa_L$.
\item
Suppose that $\fa_j$ is given.
Take a good resolution $f_j \colon X_j \to X$ for $\fa_j$
which factors through $f$.
\item
Let $S_j+B_j$ denote the birational transform of $S+B$ on $X_j$,
and $\{E_{j,i}\}_{i \in I_j}$ the set of prime divisors on $X_j$
which appears in either the exceptional locus of $f_j$ or the support of $B_j$.
We can regard the set $I$ of indices as a subset of $I_j$
by setting $E_{j,i}$ for $i \in I$
as the birational transform of $E_i$.
Set $I_{j,S} := \{i \in I_j \mid E_{j,i} \cap S_j \neq \emptyset\}$.
\item
Set $a_i:=a_{K_X+S+B}(v_{E_{j,i}})$ for $i \in I_j$.
Then $K_{X_j}+S_j = f_j^*(K_X+S+B)+\sum_{i \in I_j}a_iE_{j,i}$.
Write $\fa_j\cO_{X_j}=\cO_{X_j}(-\sum_{i \in I_j}n_{j,i}E_{j,i})$.
We set
$\fa_{j+1}:=
\fa_L \cap f_{j*}\cO_{X_j}(\sum_{i \in I_j}(\rd{a_i}+1-n_{j,i})E_{j,i})$.
Here $\rd{a}$ denotes the \textit{round-down} of $a$,
that is, the largest integer which is at most $a$.
\end{enumerate}

The below inequality is obvious by the construction of $\fa_{j+1}$.

\begin{lemma}\label{lem:n}
$n_{j+1,i} \ge n_{j,i}-(\rd{a_i}+1)$ for $i \in I$.
\end{lemma}

Let $\fa_S$ denote the ideal sheaf on $X$ defining $S$.
Kawamata--Viehweg vanishing theorem \cite[Theorem 1-2-3]{KMM87} deduces
the following equality.

\begin{lemma}\label{lem:S}
$\fa_j+\fa_S=\fa_L$ for any $j$.
\end{lemma}

\begin{proof}
We shall prove the lemma by induction on $j$.
$\fa_0+\fa_S=\fa_L$ is trivial by $\fa_0=\fa_L \supset \fa_S$.
Suppose that $\fa_j+\fa_S=\fa_L$ for fixed $j$.
Take a small positive rational number $\epsilon$
so that $\ru{a_i-(1-\epsilon)n_{j,i}}=\rd{a_i}+1-n_{j,i}$
for any $i \in I_j$ with $n_{j,i}>0$,
where $\ru{a}$ denotes the \textit{round-up} of $a$,
that is, the smallest integer which is at least $a$.
We have a natural exact sequence
\begin{multline*}
0 \to \cO_{X_j}(\sum_{i \in I_j}\ru{a_i-(1-\epsilon)n_{j,i}}E_{j,i}-S_j) \to \\
\cO_{X_j}(\sum_{i \in I_j}\ru{a_i-(1-\epsilon)n_{j,i}}E_{j,i}) \to
\cO_{S_j}(\sum_{i \in I_{j,S}}\ru{a_i-(1-\epsilon)n_{j,i}}E_{j,i}|_{S_j}) \to
0.
\end{multline*}
Since
$\sum_{i \in I_j}(a_i-(1-\epsilon)n_{j,i})E_{j,i}-S_j \equiv_{f_j}
K_{X_j}-(1-\epsilon)\sum_{i \in I_j}n_{j,i}E_{j,i}$,
we have
$R^1f_{j*}\cO_{X_j}(\sum_{i \in I_j}\ru{a_i-(1-\epsilon)n_{j,i}}E_{j,i}-S_j)=0$
by Kawamata--Viehweg vanishing theorem \cite[Theorem 1-2-3]{KMM87}.
Here we should mention that this theorem \cite[Theorem 1-2-3]{KMM87} holds
even if $D$ is an $\bR$-divisor by its proof.
Hence we obtain the surjective map
\begin{align*}
f_{j*}\cO_{X_j}(\sum_{i \in I_j}\ru{a_i-(1-\epsilon)n_{j,i}}E_{j,i}) \to
f_{j*}\cO_{S_j}
(\sum_{i \in I_{j,S}}\ru{a_i-(1-\epsilon)n_{j,i}}E_{j,i}|_{S_j}).
\end{align*}
For any $i \in I_{j,S}$,
we have $\ru{a_i-(1-\epsilon)n_{j,i}} \ge -n_{j,i}$
unless $a_i=-1$ with $n_{j,i}=0$,
and if $a_i=-1$ then $f_j(E_{j,i}|_{S_j}) \subset L$.
This implies that
\begin{align*}
f_{j*}\cO_{S_j}(\sum_{i \in I_{j,S}}\ru{a_i-(1-\epsilon)n_{j,i}}E_{j,i}|_{S_j})
&\supset f_{j*}\cO_{S_j}(-\sum_{i \in I_{j,S}}n_{j,i}E_{j,i}|_{S_j})
\cap \fa_L\cO_S \\
&\supset \fa_j\cO_S \cap \fa_L\cO_S = \fa_L\cO_S.
\end{align*}
Here we used the assumption that $\fa_j+\fa_S=\fa_L$.
On the other hand
$\ru{a_i-(1-\epsilon)n_{j,i}} \le \rd{a_i}+1-n_{j,i}$ for any $i \in I_j$,
whence
\begin{align*}
f_{j*}\cO_{X_j}(\sum_{i \in I_j}\ru{a_i-(1-\epsilon)n_{j,i}}E_{j,i})
\subset f_{j*}\cO_{X_j}(\sum_{i \in I_j}(\rd{a_i}+1-n_{j,i})E_{j,i})
\subset \cO_X.
\end{align*}
Therefore the above surjective map induces another surjective map
\begin{align*}
\fa_{j+1}=
\fa_L \cap f_{j*}\cO_{X_j}(\sum_{i \in I_j}(\rd{a_i}+1-n_{j,i})E_{j,i})
\to \fa_L\cO_S.
\end{align*}
Hence $\fa_{j+1}+\fa_S=\fa_L$ and Lemma \ref{lem:S} is proved.
\end{proof}

Let $\fm$ denote the ideal sheaf on $X$ defining $\{x\}$.
The following lemma is due to Zariski's subspace theorem \cite[(10.6)]{A98}.

\begin{lemma}\label{lem:ZST}
Unless $a_i\ge-1$ for all $i \in I$,
we can find an integer $j(l)$ with $\fa_{j(l)} \subset \fm^l$ for each $l$.
\end{lemma}

\begin{proof}
We have $\fa_j \subset f_*\cO_{\bar{X}}(-n_{j,i}E_i)$ for any $i \in I$
by the construction of $\fa_j$.
Since $n_{j,i} \ge n_{0,i}-j(\rd{a_i}+1)$ by Lemma \ref{lem:n},
we obtain that $\fa_j \subset f_*\cO_{\bar{X}}(-(n_{0,i}-j(\rd{a_i}+1))E_i)$.
Hence Lemma \ref{lem:ZST} is reduced to finding an integer $k(l)$
with $f_*\cO_{\bar{X}}(-k(l)E_i) \subset \fm^l$ for fixed $i$ and $l$.
Take a closed point $y$ in $E_i \cap f^{-1}(x)$.
Then the existence of $k(l)$ is obtained by
Zariski's subspace theorem \cite[(10.6)]{A98}
for $R:=\cO_{X,x}$ and $S:=\cO_{\bar{X},y}$.
Note that $\cO_{X,x}$ is analytically normal;
see \cite[(37.5)]{N62} for the proof.
\end{proof}

We shall complete the proof of the theorem.
Unless $a_i\ge-1$ for all $i \in I$,
we have $\bigcap_{j\ge0}(\fa_j+\fa_S)=\fa_S$
by Lemma \ref{lem:ZST} and the Artin--Rees lemma;
so $\fa_L=\fa_S$ by Lemma \ref{lem:S},
which is absurd.
Therefore $a_i\ge-1$ for all $i \in I$ and the theorem is proved.

\bibliographystyle{amsplain}

\begin{thebibliography}{9}
\bibitem{A98}
S. Abhyankar,
\textit{Resolution of singularities of embedded algebraic surfaces},
Second edition, Springer Monographs in Mathematics, Springer-Verlag (1998).
\bibitem{H64}
H. Hironaka,
Resolution of singularities of an algebraic variety
over a field of characteristic zero I, II,
Ann.\ of Math.\ (2) \textbf{79} (1964), 109-203; ibid.\ 205-326.
\bibitem{KMM87}
Y. Kawamata, K. Matsuda and K. Matsuki,
Introduction to the minimal model problem,
\textit{Algebraic geometry},
Adv.\ Stud.\ Pure Math.\ \textbf{10} (1987), 283-360.
\bibitem{K+92}
J. Koll\'ar et al,
\textit{Flips and abundance for algebraic threefolds},
Ast\'erisque \textbf{211} (1992).
\bibitem{KM98}
J. Koll\'ar and S. Mori,
\textit{Birational geometry of algebraic varieties},
Cambridge Tracts in Mathematics \textbf{134},
Cambridge University Press (1998).
\bibitem{N62}
M. Nagata,
\textit{Local rings},
Interscience Tracts in Pure and Applied Mathematics \textbf{13},
Interscience Publishers a division of John Wiley \& Sons (1962).
\bibitem{S92}
V. Shokurov,
Three-dimensional log perestroikas,
Izv.\ Ross.\ Akad.\ Nauk Ser.\ Mat.\ \textbf{56} (1992), 105-203;
translation in Russian Acad.\ Sci.\ Izv.\ Math.\ \textbf{40} (1993), 95-202.
\end{thebibliography}

\end{document}